1

# THE SCHAPER FORMULA AND THE LASCOUX, LECLERC AND THIBON-ALGORITHM

STEEN RYOM-HANSEN

ABSTRACT. We deduce the Schaper formula for Hecke-algebras at root of unity from the Jantzen conjecture in the LLT-setup. This explains an observation due to R. Rouquier.

## 1. Introduction

The Schaper formula for the determinant of the Gram matrix is one of the most useful tools for calculating decomposition numbers of the symmetric group in positive characteristic. It was first obtained by James-Murphy [JM] within the framework of the symmetric group and later derived by Schaper from the Jantzen sum formula for algebraic groups in type A, [S]. He also related it to a sum formula for a filtration of the Specht module, along the lines of Jantzen's original work.

As a matter of fact the two formulas are completely equivalent. The difference in appearance of the two formulas becomes visible first of all in case of small primes (or very singular weights) where the Jantzen sum formula has many terms that cancel out. Thus in all such situations, the Schaper formula is the more effective tool. For instance, the Carter conjecture, which describes Specht modules that remain simple under reduction mod $p$, needed the Schaper version of the formula to be solved.

The setup of this note will be the representation theory of Hecke algebras at root of unity. James and Mathas [JMa] have generalized the methods of [JM] to obtain a Schaper formula for the Specht modules in this situation. Our aim is to use the Fock module approach to representation theory in type A, due to Lascoux, Leclerc

[1]Supported by EPSRC grant M22536 and by the TMR-network algebraic Lie Theory ERB-FMRX-CT97-0100

and Thibon [LLT] to explain why the formula looks the way it does, focusing on the hooks lengths of the formula. We assume the Jantzen conjecture in the situation and deduce quickly the Schaper formula. This explains an observation of R. Rouquier presented in [LLT], which, as a matter of fact, was a starting point of our work.

## 2. **Preliminaries**

In this section we setup the notation around the LLT algorithm and state the Schaper sum formula. The $q$–analogue of the Fock space is

$$\mathcal{F}_q = \bigoplus_{\lambda \in \mathrm{Par}} \mathbb{Q}(q)|\lambda\rangle$$

with basis $|\lambda\rangle$ parametrized by the set of all partitions Par. It was originally introduced by Hayashi in [H]. Following [LT] we can interprete the basis element $|\lambda\rangle$ as the semiinfinite wedge $u_I = u_{i_1} \wedge u_{i_2} \wedge \dots$ where the sequence $I = (i_1, i_2, \dots)$ is constructed from $\lambda = (\lambda_1, \lambda_2, \dots, \lambda_n)$ using the rule $i_k = \lambda_k - k + 1$. The basis elements of $\mathcal{F}_q$ in this way correspond to the set $\mathcal{J}$ of decreasing sequences $I = (i_1, i_2, \dots)$ such that $i_k = -k + 1$ for $k$ large enough.

In this setting Leclerc and Thibon [LT] define (modulo a couple of misprints) an involution of the Fock space in the following way: For $m \geq 0$ denote by $\mathcal{J}_m$ the subset of $\mathcal{J}$ consisting of those $I$ such that $\sum_k (i_k + k - 1) = m$. For $I \in \mathcal{J}_m$ denote by $\alpha_{n,k}(I)$ the number of pairs $(r, s)$ with $1 \leq r < s \leq k$ and $r - s \not\equiv 0 \bmod n$. Then the bar involution of the Fock space is defined by $\overline{q} = q^{-1}$ and

$$\overline{u_I} = (-1)^{\binom{k}{2}} q^{\alpha_{n,k}(I)} u_{i_k} \wedge u_{i_{k-1}} \wedge \cdots \wedge u_{i_1} \wedge u_{i_{k+1}} \wedge u_{i_{k+2}} \wedge \cdots \tag{1}$$

for (any) $k \geq m$.

The straightening rules to bring this in normal form (with decreasing indices) are the following: Suppose that $l < m$ and that $m - l = i \bmod n$. If $i = 0$

$$u_l \wedge u_m = -u_m \wedge u_l \tag{2}$$

and otherwise

$$u_l \wedge u_m = -q^{-1} u_m \wedge u_l + (q^{-2}-1)\big(u_{m-i} \wedge u_{l+i} \\ -q^{-1} u_{m-n} \wedge u_{l+n} + q^{-2} u_{m-n-i} \wedge u_{l+n+i} - \cdots\big) \tag{3}$$

Let $\mu \vdash m$. Define $\alpha_{\lambda,\mu}(q)$ by

$$\overline{|\mu\rangle} = \sum_{\lambda \vdash m} \alpha_{\lambda,\mu}(q) |\lambda\rangle. \tag{4}$$

Then by [LT] we have $\alpha_{\lambda,\mu} \in \mathbb{Z}[q,q^{-1}]$, $\alpha_{\lambda,\mu} = 0$ unless $\lambda \trianglelefteq \mu$ and $\alpha_{\lambda,\lambda} = 1$.

The Fock space can be made into an integrable module for the quantum group $U_q(\widehat{\mathfrak{gl}_n})$ and thus possesses a lower crystal basis, see [LT] for an explanation of this. Let $L$ be the $\mathbb{Z}[q]$-sublattice of $\mathcal{F}_q$ with basis $\{|\lambda\rangle \mid \lambda \in \text{Par}\}$. The lower global basis $G(\lambda)$ of $\mathcal{F}_q$ can now be defined by the following conditions

$$\overline{G(\lambda)} = G(\lambda), \quad G(\lambda) = |\lambda\rangle \bmod qL \tag{5}$$

Let $d_{\mu\lambda}(q)$ be defined by

$$G(\lambda) = \sum_{\mu} d_{\mu\lambda}(q) |\mu\rangle \tag{6}$$

Then $d_{\mu\lambda}(q) \in \mathbb{Z}[q]$, $d_{\mu,\lambda}(q) = 0$ unless $\lambda \trianglelefteq \mu$ and $d_{\lambda,\lambda}(q) = 1$.

Let $S(\lambda)$ be the Specht module for the Hecke algebra of type A on $m$ letters, specialized at an $n$'th root of unity and let for $\lambda$ a $n$-regular partition of $m$, $D(\lambda)$ be its unique simple quotient, see eg. [LLT] for an definition of these concepts.

As explained in for example [JMa], there is a natural bilinear form on the integral Specht module $S_R(\lambda)$ which gives rise to a Jantzen filtration and a sum formula on $S(\lambda)$, provided $R$ is a principal ideal domain with a prime $\mathfrak{p}$. We take $R$ as in (4.2)(i) in loc. cite., i.e. $R := \mathbb{F}[q,q^{-1}]$, for $\mathbb{F}$ any field and $\mathfrak{p}$ a cyclotomic polynomial in $R$ such that $q$ has order $n$ in the residue field $R/\mathfrak{p}$. Denoting by

$$S(\lambda) = S(\lambda)^0 \supseteq S(\lambda)^1 \supseteq S(\lambda)^2 \supseteq \cdots$$

the Jantzen filtration of $S(\lambda)$, LLT conjectured – inspired by an observation due to R. Rouquier – that

$$d_{\lambda\mu}(q) = \sum_{i>0} [S(\lambda)^i/S(\lambda)^{i+1} : D(\mu)]\, q^i \tag{7}$$

This conjecture corresponds to the Jantzen conjecture for Verma modules and is reported proved by Grojnowski.

Now recall that one associates to each partition $\lambda = (\lambda_1, \cdots, \lambda_r)$ a set of $\beta$-numbers generalizing the first-column hook lengths $h_1 > h_2 > \cdots > h_r$ of $\lambda$. Formally they are defined the following way: if $\{\beta_1, \beta_2, \ldots, \beta_s\}$ is a set of $\beta$-numbers for $\lambda$ then also $\{\beta_1+1, \beta_2+1, \ldots, \beta_s+1, 0\}$ is a set of $\beta$-numbers for $\lambda$ and $\{h_1, h_2, \ldots, h_r\}$ is a set of $\beta$-numbers for $\lambda$. One gets $\lambda$ back from any set of $\beta$-numbers $\{\beta_1, \beta_2, \cdots, \beta_s\}$ for $\lambda$ by choosing $\sigma \in \Sigma_s$ such that $\beta_{\sigma(1)} > \beta_{\sigma(2)} > \ldots > \beta_{\sigma(s)}$ and setting $\lambda_i = \beta_{\sigma(i)} + i - s$. We then introduce the numbers $d(\beta_1, \beta_2, \ldots, \beta_s)$:

$$d(\beta_1, \beta_2, \ldots, \beta_s) := (-1)^\sigma \dim S(\lambda)$$

Letting $\det \lambda$ denote the determinant of the bilinear form on $S_R(\lambda)$ and $\nu_{\mathfrak{p}}(\cdot)$ the $\mathfrak{p}$-adic valuation of $R$, the Schaper formula states that

$$\nu_{\mathfrak{p}}(\det \lambda) = \sum_{1\le a\le b\le s} \sum_{c=1}^{\lambda_b} \nu_{\mathfrak{p}}\left(\frac{[h_{ac}]}{[h_{bc}]}\right) d(h_{11}, h_{21}, \ldots, h_{a1}+h_{bc}, \ldots, h_{b1}-h_{bc}, \ldots, h_{s1})$$

Jantzen's classical argument then gives the following sum formula for the Jantzen filtration:

$$\sum_i S(\lambda)^i = \sum_{1\le a\le b\le s} \sum_{c=1}^{\lambda_b} \nu_{\mathfrak{p}}\left(\frac{h_{ac}}{h_{bc}}\right) S(h_{11}, h_{21}, \ldots, h_{a1}+h_{bc}, \ldots, h_{b1}-h_{bc}, \ldots, h_{s1}) \tag{8}$$

(in the Grothendieck group of Hecke modules) where $S(h_{11}, h_{21}, \ldots, )$ is the virtual Specht module obtained from the $\beta$-numbers in the analogous way.

We illustrate the hook length $h_{bc}$, the length of the hook defined by the $(b,c)'$th entry in the Young diagram, in the following example

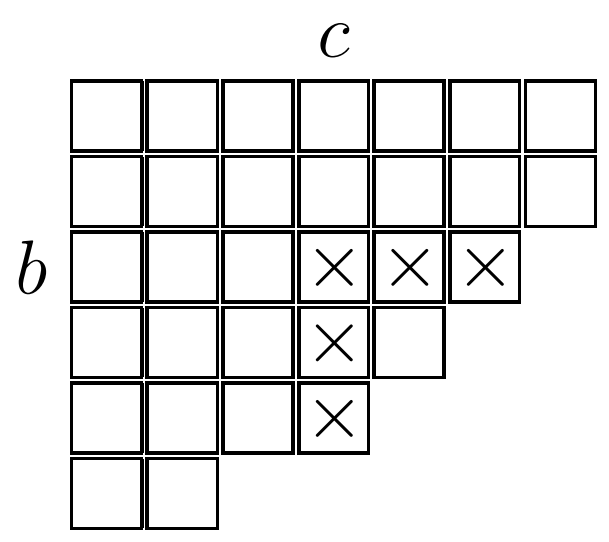

## 3. Gabber-Joseph's approach to the sum formula

The following arguments are parallel to ideas in the work of Gabber-Joseph [GJ]. Combining the formulas (4), (5) and (6) we arrive at the equality

$$d_{\lambda\mu}(q) = \sum_{\lambda \leq \tau \leq \mu} a_{\lambda\tau}(q) d_{\tau\mu}(q^{-1}) \tag{9}$$

Differentiating this we obtain the following expression

$$d'_{\lambda\mu}(1) = \frac{1}{2} \sum_{\lambda \leq \tau \leq \mu} a'_{\lambda\tau}(1)\, d_{\tau\mu}(1) \tag{10}$$

Here we used the facts that $a_{\lambda\tau}(1) = 0$ unless $\lambda = \tau$, and $a_{\lambda\lambda}(q) = 1$ which are immediate from the definition. Now, Ariki [A] has proved the following deep theorem:

$$d_{\tau\mu}(1) = [S(\tau), D(\mu)] \qquad \mu \text{ regular} \tag{11}$$

Combining this with (10) and summing over all regular $\mu$, we obtain the following expression in the Grothendieck group

$$\sum_{\mu} d'_{\lambda\mu}(1)\, [D(\mu)] = \frac{1}{2} \sum_{\lambda \leq \tau} a'_{\lambda\tau}(1)\, [S(\tau)] \tag{12}$$

Assuming the Jantzen conjecture (7) we find that

$$d'_{\lambda\mu}(1) = \sum_{i>0} i \left[ S(\lambda)^i / S(\lambda)^{i+1}, D(\mu) \right] \tag{13}$$

Once again we sum over all $\mu$ and find that

$$\sum_{\mu} d'_{\lambda\mu}(1)[D(\mu)] = \sum_{i>0} i\left[S(\lambda)^i/S(\lambda)^{i+1}\right] \tag{14}$$

The sum appearing in this expression is the left hand side of the sum formula. The idea is now to calculate $a'_{\lambda\tau}(1)$ and then (12) becomes the right hand side of the sum formula.

## 4. Calculation of $\mathbf{a}'_{\lambda\tau}(\mathbf{1})$

Although the $a_{\lambda\tau}(q)$ are combinatorially difficult to describe, it is at least clear from the definition and the straightening rules that they take the form $\pm q^{-n}(q^{-2}-1)^m$ for some $m, n \in \mathbb{N}$.

Now $m = 0$ corresponds to the case $\lambda = \tau$; but in that case we have that $a_{\lambda\tau}(q) = 1$ so $a'_{\lambda\tau}(1) = 0$.

The case $m \geq 2$ is also easily handled. We get once again that $a'_{\lambda\tau}(1) = 0$, since

$$\left((q^{-2}-1)^2 g(q^{-1})\right)'(1) = 0$$

for any polynomial $g$.

We then finally consider the case $m = 1$. When we apply the straightening rules (2) and (3), these are the terms that arise when we choose exactly once a $q$-multiple of $(q^{-2}-1)$. For a given $\tau$, the relevant $\lambda' s$ are those whose $q$-wedges are equal to those of $\tau$ except for exactly two factors $u_l, u_m$, that have been changed into $u_{l-i}, u_{m+i}$ or $u_{l-n}, u_{m+n}$ and so on.

For any $\mu \in \text{Par}$, denote by $u_\mu$ the semiinfinite wedge corresponding to the partition $\mu$.

Notice that there is a close connection between the $\beta$-numbers formalism and the wedge formalism. Indeed, if $I = (i_1, i_2, \dots) \in \mathcal{J}_m$ and $u_I = u_\tau$, for some $\tau \in \text{Par}$, one easily checks that $\{i_1 + m, i_2 + m, \dots, i_m + m, 0\}$ is a set of $\beta$-numbers for $\tau$.

Now a set of $\beta$-numbers for the partition $\lambda$ is the set of first column hook lengths of $\{h_{11}, \dots, , h_{s1}\}$ of $\lambda$. By our assumption a set of $\beta$-numbers for the partition $\tau$ is given by $\{h_{11}, \dots, h_{a1} +$ $r, \dots, h_{b1} - r, \dots, h_{s1}\}$ for some integer $r$. It then follows from the theory of $\beta$-numbers, see [JK], that $r$ equals some hook length of $\lambda$, say $r = h_{bc}$ (geometrically, $\lambda$ arises from $\tau$ by moving a skew hook of length $h_{bc}$ from $a$'th to $b$'th row).

Applying the straightening rules to $\overline{u_\tau}$ we only once do not take the commuting terms, namely for the wedge indices $m = h_{a1} + h_{bc}$ and $l = h_{b1} - h_{bc}$. From this we get two possible kinds of terms, namely those whose wedge index has been changed by a multiple of $n$ and those whose index has been changed by $i \equiv m - l \mod n$.

Now the first kind of terms correspond to the denominator of the sum formula satisfying $\nu_{\mathfrak{p}}([h_{bc}]) = 1$, while the second kind correspond to the nominator satisfying $\nu_{\mathfrak{p}}([h_{ac}]) = 1$.

Let us be more precise. The statement is clear for the first ones: passsing from $u_\tau$ to $u_\lambda$, $u_m = u_{h_{a1}+h_{bc}}$ in $u_\tau$ becomes $u_{h_{a1}}$ in $u_\lambda$ and thus the difference $h_{bc}$ is a multiple of $n$ if and only if $\nu_q([h_{bc}]) = 1$.

For the second kind of terms consider the Young diagram of $\lambda$:

$Y(\lambda) =$ 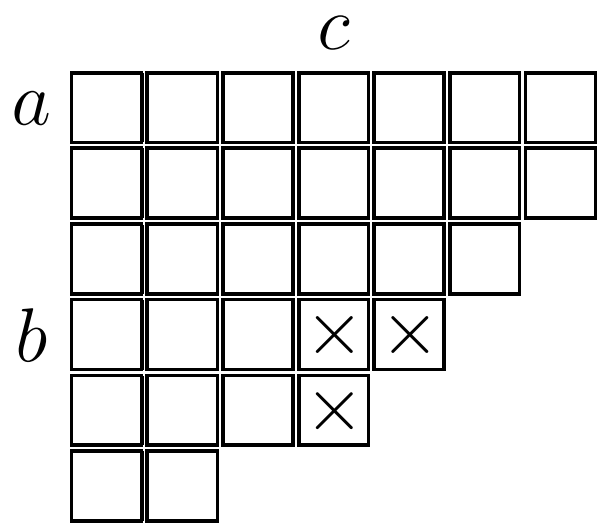

Now if $\nu_{\mathfrak{p}}([h_{ac}]) = 1$ one sees at the diagram that the difference $h_{a1} - h_{b1}$ has length $-h_{bc} \bmod n$ since it is the negative of the number of marked boxes in it. But then

$$i \equiv m - l \equiv (h_{a1} + h_{bc}) - (h_{b1} - h_{bc}) \equiv h_{bc}$$

And conversely, if $i \equiv h_{bc}$ we must have $\nu_{\mathfrak{p}}([h_{ac}]) = 1$. So indeed, the second kind of terms correspond to $\nu_{\mathfrak{p}}([h_{ac}]) = 1$.

It only remains to check that everything has been normed the right way. Here only the signs need to be checked, since the $q$-powers are of no importance for $a'_{\lambda\tau}(1)$. Now if $h_{11} > h_{21} > \cdots >$ $h_{a1} + h_{bc} > \cdots > h_{b1} - h_{bc} > \cdots > h_{s1}$ we need to perform $\binom{k}{2}$ ($k$ from the definition of the bar involution) commutations when going from $u_\tau$ to $u_\lambda$. Each of these commutations comes with a sign, which however is cancelled out with the signs from the definition of the bar involution.

If the hooks are not nicely ordered as above we use some $\sigma \in \Sigma_k$ to order them, this changes the sign of the semiinfinite wedge by a factor $(-1)^\sigma$ (because we only use the pure commutators) and this cancels out with the sign in the definition of $d(h_{11}, h_{21}, \cdots, h_{a1} + h_{bc}, \cdots, h_{b1} - h_{bc}, \cdots, h_{s1})$.

All in all we have proved the following theorem:

**Theorem 1.** *The Jantzen conjecture*

$$d_{\lambda\mu}(q) = \sum_{i>0} [S(\lambda)^i / S(\lambda)^{i+1} : D(\mu)] q^i$$

*is compatible with the sum formula*

$$\nu_q (\det \lambda) = \sum_{1 \le a \le b \le s} \sum_{c=1}^{\lambda_b} \nu_{\mathfrak{p}} \left( \frac{[h_{ac}]}{[h_{bc}]} \right) d(h_{11}, h_{21}, \ldots, h_{a1} + h_{bc}, \ldots, h_{b1} - h_{bc}, \ldots, h_{s1})$$

*for quantum groups at a n'th root of unity.*

## References

[A] S. Ariki, On the decomposition numbers of the Hecke algebra $G(m, 1, n)$, Journal of Mathematics of Kyoto University, **36**, (1996), 789-808

[GJ] O. Gabber, A. Joseph, Towards the Kazhdan–Lusztig conjecture, Ann. scient. Éc. Norm. Sup., 4$^e$ série, t.**14**, (1981), p.261 à 302

[H] T. Hayashi, $q$-analogues of Clifford and Weyl algebras– spinor and oscillator representations of quantum enveloping algebras, Commun. Math. Phys., **127**(1990), 129-144

[K] M. Kashiwara, On crystal bases of the $q$-analogue of universal enveloping algebras, Duke Math. J. **63** (1991), 465-516

[JM] G. D. James, G.E. Murphy, The Determinant of the Gram Matrix for a Specht Module, Journal of Algebra **59** (1979) 222-235.

[JMa] G. D. James, A. Mathas, A $q$-analogue of the Jantzen–Schaper theorem, Proc. Lond. Math. Soc., **74** (1991), 241-274

[JK] G. D. James, A. Kerber, The representation theory of the symmetric group, **16**, Encyclopedia of Mathematics, Addison-Wesley, Massachusetts, (1981)

[LLT] A. Lascoux, B. Leclerc, J.-Y. Thibon, Hecke algebras at roots of unity and crystal bases of quantum affine algebras, Commun. Math. Phys. **181** (1996), 205-263

[LT] B. Leclerc, J.-Y. Thibon, Int. Math. Res. Notices, **9** (1996), 447-456

[MW] K.C. Misra, T.Miwa, Crystal base for the basic representation of $U_q(\widehat{\mathfrak{sl}_\mathfrak{n}})$, Commun. Math. Phys., **134** (1990), 79-88

[S] Schaper, Characterformel für Weyl-Moduln und Specht-Moduln in Primzahlcharacteristik, Diplomarbeit, Bonn, (1981).

Matematisk Afdeling, Universitetsparken 5, DK-2100 København Ø, Danmark